\documentclass{amsart}
\usepackage{amssymb, amscd}
\usepackage{enumerate}

\newcommand\blue[1]{{#1}}
\newcommand\datver[1]{\def\datverp
 {\par\boxed{\boxed{\text{Version: #1; Run: \today}}}}}
\datver{4.0; Revised: Nov 30, 2005 by Bernd}
\newcommand{\ie}{{\em i.\thinspace e., }}

\newcommand{\CI}{\mathcal{C}^\infty}
\newcommand{\CIc}{{\mathcal C}^{\infty}_{\text{c}}}
\newcommand{\pa}{\partial}

\newcommand{\maC}{\mathcal C}

\newcommand{\maP}{\mathcal P}

\newcommand{\CC}{\mathbb C}

\newcommand{\PP}{\mathbb P}
\newcommand{\RR}{\mathbb R}

\newcommand{\ZZ}{\mathbb Z}

\newcommand{\dvol}{\mathop{{\rm dvol}}}

\newcommand\comment[1]{}

\newcommand{\maK}{\mathcal K}

\newcommand{\DiffO}[1]{{\rm Diff}^{#1}_{0}}

\newcommand{\supp}{\operatorname{supp}}
\newcommand\Kond[2]{{\mathcal K}^{#1}_{#2}}
\newcommand{\BK}{Babu\v{s}ka--Kondratiev}
\newcommand{\CIS}{\maC^{\infty}(\Sigma\PP)}

\newtheorem{theorem}{Theorem}[section]
\newtheorem{proposition}[theorem]{Proposition}
\newtheorem{corollary}[theorem]{Corollary}

\newtheorem{lemma}[theorem]{Lemma}

\theoremstyle{definition}

\theoremstyle{remark}

\author[B. Ammann]{Bernd Ammann} \address{Bernd Ammann,
Institut \'Elie Cartan, Universit\'e Henri Poincar\'e Nancy 1, B.P. 239,
 54506 Vandoeuvre-Les-Nancy, France,
{\rm http://www.differentialgeometrie.de/ammann}}
\email{bernd.ammann@gmx.net}

\author[V. Nistor]{Victor Nistor} \address{Victor Nistor, Pennsylvania State
University, Math. Dept., University Park, PA 16802, USA}
\email{nistor@math.psu.edu}

\begin{document}
\title[Sobolev spaces]{Weighted Sobolev spaces and regularity for
polyhedral domains}

\dedicatory{\begin{center}{\em Dedicated to Ivo Babu\v{s}ka on the
occasion of his 80th birthday.}\end{center}}

\date{\today}

\begin{abstract} We prove a regularity result for
the Poisson problem $-\Delta u = f$, $u \vert_{\pa \PP} = g$ on a
polyhedral domain $\PP \subset \RR^3$ using the \BK\ spaces
$\Kond{m}{a}(\PP)$. These are weighted Sobolev spaces in which the
weight is given by the distance to the set of edges \cite{Babu70,
Kondratiev67}. In particular, we show that there is no loss of
$\Kond{m}{a}$--regularity for solutions of strongly elliptic
systems with smooth coefficients. We also establish a ``trace
theorem'' for the restriction to the boundary of the functions in
$\Kond{m}{a}(\PP)$.
\end{abstract}

\maketitle

\tableofcontents

\section{Introduction}

Let $\Omega \subset \RR^n$ be a smooth, bounded domain. Then it is
well known \cite{BabuAziz, BrennerScott, Evans, Necas, Taylor1}
that the equation
\begin{equation}\label{eq.BVP}
    \Delta u = f \in H^{m-1}(\Omega), \quad u = 0\; \text{ on }\; \pa
    \Omega,
\end{equation}
has a unique solution $u \in H^{m+1}(\Omega)$. In particular, $u$
will be smooth on $\overline{\Omega}$ if $f$ is smooth on
$\overline{\Omega}$. This well-posedness result is especially
useful in practice for the numerical approximation of the solution
$u$ of Equation \eqref{eq.BVP}, see for example \cite{BabuAziz,
BabuSobolev, BrennerScott} among many possible references.

In practice, however, it is rarely the case that $\Omega$ is
smooth. In fact, if $\pa \Omega$ is {\em not} smooth, then the
smoothness of~$f$ on $\overline{\Omega}$ does not imply that the
solution~$u$ of Equation \eqref{eq.BVP} is also smooth on
$\overline{\Omega}$. Therefore there is a {\em loss of regularity} for
elliptic problems on non-smooth domains. Wahlbin \cite{Wahlbin84} (see
also \cite{Babu71, Laasonen57, Wahlbin91}) has shown that this leads
to some inconvenience in numerical applications, namely that a
quasi-uniform sequence of triangulations on $\Omega$ will {\em not}
lead to optimal rates of convergence for the Galerkin approximations
$u_h$ of the solution of \eqref{eq.BVP}.

The loss of regularity can be avoided, however, if one removes the
singular points by ``sending them to infinity'' by suitably changing
the metric with a conformal factor. It can be proved then that the
resulting Sobolev spaces are the ``Sobolev spaces with weights''
considered for instance in \cite{BabuAziz, BabuRo, BabuSobolev,
Kondratiev67} \blue{and in several other papers. A related
construction, leading however to countably normed spaces}, was
considered in \cite{BabuGuo}. Let $f> 0$ be a smooth function on a
domain $\Omega$, then the \emph{$m$th Sobolev space with weight $f$}
is defined by
\begin{equation}\label{eq.def.WS}
    \maK^{m}_{a}(\Omega; f) := \{ u, \
    f^{|\alpha|-a} \pa^{\alpha} u \in L^2(\Omega), \
    |\alpha| \le m \}, \qquad m \in \ZZ_+,\ a \in \RR.
\end{equation}
The regularity result for Equation \eqref{eq.BVP} extends to
polyhedral domains~$\PP$ in three dimensions with the usual
Sobolev spaces replaced by the spaces $\maK^{m}_{a}(\PP) :=
\maK^{m}_{a}(\PP; \vartheta)$, $\vartheta$ being the distance to
the edges. The spaces $\maK^{m}_{a}(\pa \PP; \vartheta)$ on the
boundary are defined similarly for $m \in \ZZ_+ := \{0, 1, \ldots
\}$; for $m \in \RR_+$ they are defined using interpolation.

\begin{theorem}\label{thm.princ}
Let $\PP \subset \RR^3$ be a polyhedral domain. Let $m \in \ZZ_+$
and $a \in \RR$. Assume that $u \in \Kond{1}{a+1}(\PP)$, $\Delta u
\in \Kond{m-1}{a-1}(\PP)$, and $u\vert_{\pa \PP} \in
\Kond{m+1/2}{a+1/2}(\pa \PP, \vartheta)$, then $u \in
\Kond{m+1}{a+1}(\PP)$ and there exists $C>0$ independent of $u$
such that
\begin{equation*}
    \|u\|_{\Kond{m+1}{a+1}(\PP)} \le
    C\big(\|\Delta u\|_{\Kond{m-1}{a-1}(\PP)} +
    \|u\|_{\Kond{0}{a+1}(\PP)} +
    \|u \vert_{\pa \PP}\|_{\Kond{m+1/2}{a+1/2}(\pa \PP, \vartheta)} \big).
\end{equation*}
The same result holds if we replace $\Delta$ with a strongly
elliptic operator or system.
\end{theorem}

Theorem \ref{thm.princ} is well known in two dimensions, \ie for
polygonal domains, \blue{and for domains with conical points}
\cite{Camus, Kondratiev67}. \blue{See also \cite{CDMaxwell,
KMRossmann, MazyaBoth} where similar results were proved using a
dyadic partition of unity technique.} For the result in two
dimensions, $\vartheta$ is the distance to the vertices of the
polygonal domain considered or to the conical points. In general, in
$d$ dimensions, one takes $\vartheta(x)$ to be the distance to the set
of non-smooth boundary points of $\PP$.  \blue{Significantly less
papers have dealt with the case of three dimensions. Nevertheless, let
us mention the following. A general and far reaching theory (valid
also in higher dimensions) was developed by Dauge in
\cite{Dauge}. Regularity estimates based on singular function
expansions were proved by Apel and Nicaise \cite{NicaiseApel} and
Lubuma and Nicaise \cite{LubumaNicaise95}.  These results were then
applied in these papers in order to obtain optimal rates of
convergence in the Finite Element Method.  In \cite{MazyaRossmann},
Mazya and Rossmann have obtained similar results using estimates on
Green functions. Buffa, Costabel, and Dauge \cite{BCD} have proved or
stated similar regularity and well posedness results for polyhedral
dimensions in three dimensions. Our modified weight $r_\Omega$ was
introduced in \cite{CDMaxwell2}, where the above regularity theorem
was proved for $m = 1$. In \cite{KelloggOsborn}, Kellogg and Osborn
have obtained regularity results of a similar kind for the Stokes
operator. \blue{Borsuk and Kondratiev established many regularity
results for Dini-Liapunov regions in $\RR^n$, $n\geq 3$, in their
recent monograph \cite{borsuk.kondratiev}. Note that the notion of a
Dini-Liaponov region is a generalisation of a domain with
$C^{1,\alpha}$-bondary}.  See also \cite{ArnoldFalk, CostabelDauge,
Dauge2, KelloggTr1, MazyaRossmann2, NicaiseB, NicaiseM3AS}, to mention
just a few other papers. A regularity result valid in all dimensions
was obtained in \cite{AIN} using ``Lie manifolds.''}

\blue{We are grateful to one of the referees, who pointed out to
us that Theorem \ref{thm.princ} can also be obtained from the
results of the monograph \cite{NP}.} In this paper, we follow
\cite{AIN}, but we use more elementary methods that lead to a
short proof. We also introduce some ideas that are specific to
polyhedral domains in three dimensions and may be useful in
applications to Numerical Analysis. Moreover, our paper is
self-contained and the references to \cite{AIN} are only for
comparison.

We would like to stress that Theorem \ref{thm.princ} does {\em not}
constitute a Fredholm (or ``normal solvability'') result, because the
inclusion $\Kond{m+1}{a+1}(\PP) \to \Kond{0}{a+1}(\PP)$ is {\em not
compact} for all $m$ and $a$ \cite{AIN}. By contrast, if $\PP$ is a
polygon, then $P = - \Delta$ with Dirichlet boundary conditions is
a Fredholm operator from $\Kond{m+1}{a+1}(\PP)$ to
$\Kond{m-1}{a-1}(\PP)$ precisely when $a$ is different from $k
\pi/\alpha$, where $k \in \ZZ$, $k \neq 0$, and $\alpha$ ranges
through the angles of the polygon \cite{Kondratiev67, KMRossmann}.

The Poincar\'e inequality $\|u\|_{\Kond{1}{1}(\PP)} \le C
\|\nabla u\|_{L^2(\Omega)}$ proved in \cite{BNZ-I}, gives that
$\Delta$ is coercive on the space $\Kond{1}{1}(\PP)$ and hence the
map $\Delta : \Kond{1}{1}(\PP) \cap \{u = 0 \text{ on } \pa \PP \}
\to \Kond{-1}{-1}(\PP)$ is a continuous bijection. By combining
this with Theorem \ref{thm.princ} we obtain that
\begin{equation}
    \Delta : \Kond{m+1}{a+1}(\PP) \cap \{u = 0 \text{ on } \pa \PP \}
    \to \Kond{m-1}{a-1}(\PP)
\end{equation}
is a continuous bijection, for any $m \in \ZZ_+$ and $|a| < \eta$,
with $\eta$ depending only on $\PP$. The same result holds if $\Delta$
is replaced with $P + c_P$, where $P$ is a strongly elliptic system
with smooth coefficients and $c_P > 0$ and $\eta > 0$ are constants
depending \blue{only} on $P$ \cite{BNZ-I}.

To prove Theorem \ref{thm.princ}, we first introduce the weighted
Sobolev spaces $\maK^{m}_{a}(\pa \PP, \vartheta)$ on the boundary of
$\PP$. For $m \not \in \ZZ_+$, these spaces are defined by duality and
interpolation. Then we provide an alternative definition of the spaces
$\maK^{m}_{a}(\PP) := \maK^{m}_{a}(\PP, \vartheta)$ and
$\maK^{m}_{a}(\pa \PP, \vartheta)$ using partitions of unity. This
allows us to define a trace map $\maK^{m}_{a}(\PP) \to
\maK^{m-1/2}_{a-1/2}(\pa \PP, \vartheta)$, which extends the
restriction map and is a continuous surjection, as in the case of a
smooth domain. We also show that any differential operator $P$ of
order $m$ with smooth coefficients induces a continuous map $P :
\maK^{s}_{a}(\PP) \to \maK^{s-m}_{a-m}(\PP)$.

We need to introduce an enhanced space of smooth, bounded functions
$\CI(\Sigma \PP)$, which contains the cylindrical and spherical
coordinates functions and is minimal with this property.  In
particular, $\CI(\overline{\Omega}) \subset \CIS \subset
\CI(\Omega)$. Let $\rho_P(p)$ be the distance from $p$ to the vertex
$P$ of $\PP$ and $r_e(p)$ be the distance from $p$ to the line
determined by the edge $e$ of $\PP$ (for $\PP$ non-convex we need to
slightly change the definition of $r_e$). Then $\rho_P, \rho_e
\in \CIS$, although they are not smooth functions on
$\overline{\Omega}$ in the usual sense. Let $A$ and $B$ be the end
vertices of the edge $e$ (\ie $e = [AB]$). We further define $\tilde
r_{e} := \rho_A^{-1} \rho_B^{-1} r_e$ and $r_{\PP} = \prod_e \tilde
r_e \times \prod_P \rho_P$. Then $\tilde r_e, r_{\PP} \in \CIS$. The
functions in $\CIS$ have the following strong boundedness property
\begin{equation}
    (r_{\PP}\pa_x)^i (r_{\PP}\pa_y)^j (r_{\PP}\pa_z)^k u
    \in \CIS \subset L^\infty(\PP)
\end{equation}
for all $u \in \CIS$. The consideration of $\CIS$ and of the
derivatives of the form $r_{\PP}\pa_x$, $r_{\PP}\pa_y$, and
$r_{\PP}\pa_z$ is a substitute for the results on Lie manifolds
used in \cite{AIN}. However, the results of \cite{AIN} also apply
to non-compact manifolds and to a larger class of singular
domains.

The methods of this paper are used for a general regularity
and well-posedness result for anisotropic elasticity in general
polyhedral domains (including cracks) in \cite{MazzucatoNistor2}.
We do not include in this paper any concrete applications, but let
us refer the reader to \cite{NicaiseApel, ArnoldFalk, BKP, BNZ-II,
BCD, CDMaxwell2}, where concrete applications of results similar
to ours were provided.

\subsubsection*{Acknowledgements:}\ We would like to thank Ivo
Babu\v{s}ka, Constantin Bacuta, Alexandru Ionescu, Robert Lauter,
Anna Mazzucato, and Ludmil Zikatanov for useful discussions. The
first named author wants to thank MSRI, Berkeley, CA, USA for its
hospitality, the last named author thanks Institut \'Elie Cartan,
Nancy, France, where part of the work has been completed.

%%%%%%%%%%%%%%%%%%%%%%%%%%%%%%%%%%%%%%%%%%%%%%%%%%%%%%%%%%%%%%%%%%%
%%%%%%%%%%%%%%%%%%%%%%%%%%%%%%%%%%%%%%%%%%%%%%%%%%%%%%%%%%%%%%%%%%%

\section{Smooth functions and differential operators on
$\PP$\label{sec2}}

In this section, we shall introduce the space $\CIS \subset \CI(\PP)
\cap L^\infty(\PP)$ and relate it to the differentials $r_{\PP}\pa_x$,
$r_{\PP}\pa_y$, and $r_{\PP}\pa_z$ mentioned in the
Introduction. Similar vector fields have appeared also in
\cite{BNZ-ND}. \blue{When only edges are involved (\ie no vertices),
the use of these vector fields goes back to \cite{Mazzeo1,
Mazzeo2}. See also \cite{MelroseScattering, SchroheFC}.}

\subsection{Polygons and polyhedral domains}
Let us fix some terminology to be used in what follows.

A {\em polygon} $\PP_0$ in a two dimensional Euclidean space is an
open, connected subset whose boundary consists of finitely many
straight segments (possibly of infinite length) called {\em sides} and
having at most the end points in common. For simplicity, we assume
that $\pa \PP_0 = \pa \overline{\PP}_0$, \blue{which means that no
point of the boundary $\pa \PP_0$ is in the interior of
$\overline{\PP}_0$} (thus cracks are excluded). The points common to
more than one straight segment of the boundary are called the vertices
of $\PP_0$. We require that each vertex belongs to exactly two sides.

We do not require the boundary of $\PP_0$ to be connected. For
simplicity, in this paper we also assume that the sides are
maximally extended, so that they are not contained in larger
segments contained in the boundary. This assumption is however not
essential.

Similarly, a {\em polyhedral domain} $\PP \subset \RR^3$ is a
connected, open subset whose boundary satisfies $\pa \PP = \pa
\overline{\PP}= \bigcup_{j = 1}^N \overline{D}_j$ and:
\begin{enumerate}[(i)]
\item each $D_j$ is a polygon contained in an
affine $2$-dimensional subspace of $\RR^3$;
\item the sets $D_j$ are disjoint;
\item a side of $D_j$ is a side of exactly one other $D_k$.
\end{enumerate}

The vertices of $\PP$ are the vertices of the polygonal domains
$D_j$. The edges of $\PP$ are the sides of the polygonal domains
$D_j$. Hence an edge belongs to exactly two faces of $\PP$. For
each vertex $P$ of $\PP$, we choose a small open ball $V_P$
centered in $P$. We assume that the neighborhoods $V_P$ are chosen
to be disjoint.

We stress that, in our convention, both the polygons and the
polyhedra are {\em open} subsets. We do not require these sets to
be bounded in general, although this assumption is needed for some
of our results.

\subsection{Useful functions and other notation}
Assume, for the definition of $r_e$, $\theta_{e}$, and $\phi_{P, e}$
in this subsection, that $\PP$ is {\em convex}. If $\PP$ is not
convex, then we slightly change the definitions of these functions
such that the new functions retain their behaviour around $e$, but
will become smooth everywhere in space except on $\overline{e}$
\cite{AIN}. The modified functions $\phi_{P,e}$ and $\theta_e$ will
then be defined and smooth on $\PP$. We postpone the technical
construction of the modified functions $\phi_{P,e}$ and $\theta_e$ for
the Appendix, in order not to interrupt the flow of the
presentation. \blue{(Let us stress, however, that none of our results
requires the assumption that $\PP$ be convex.)}

Let us first recall from the Introduction that we have denoted by
$\rho_P(p)$ the distance from $p$ to the vertex $P$ of $\PP$.
Also, recall that we have denoted by $r_e(p)$ the distance from
$p$ to the {\em line} determined by the edge $e$ of $\PP$ and by
\begin{equation}\label{eq.rPP}
    r_{\PP} := \prod_e \tilde r_e \times \prod_P \rho_P,
    \qquad \text{where } \; \tilde r_{e} :=
    \rho_A^{-1} \rho_B^{-1} r_e\; \text{ for }\; e = [AB].
\end{equation}
In the above formula, the products are taken over all vertices $P$
and all edges $e$ of~$\PP$. The notation $e = [AB]$ means that $e$
is the edge joining the vertices $A$ and $B$. If $e = [A,
\infty)$, that is, if $e$ is a half-line, then $\tilde r_{e} :=
\rho_A^{-1} r_e$. Finally, if $e$ is infinite in both directions
(\ie for a dihedral angle), we let $\tilde r_{e} := r_e$.

Choose for each edge $e$ a plane $\maP_e$ containing one of the faces
$D_j$ of $\PP$ such that $e \subset \overline{D}_j$. If $x$ is not on
the line defined by $e$, we define $\theta_e$ to be the angle in a
cylindrical coordinates system $(r_e, \theta_e, z)$ determined by the
edge $e$ and the plane $\maP_e$. More precisely, let $q \in e$ be the
foot of the perpendicular from $p$ to $e$.  Then $\theta_e(p)$ is the
angle between $pq$ and $\maP_e$.  Similarly, for each vertex $P$ and
edge $e$ adjacent to $P$, we define $\phi_{P, e}(p)$ to be the angle
between the segment $pP$ and the edge $e$ (except for $p=P$, in which
case $\phi_{P, e}(p)$ is not defined).

If $\PP$ is convex, then the functions $\theta_e$ and $\phi_{P,e}$
are defined and smooth on $\PP$ (recall that $\PP$ is an open
subset). \blue{They will be part of the spherical coordinate
system $(\rho_P, \theta_e, \phi_{P, e})$ centered at $P$.} For
$\PP$ non-convex, this property will be enjoyed by the modified
functions $\theta_e$ and $\phi_{P, e}$ introduced in the Appendix.
All the following definitions and constructions below are the same
in the case of a non-convex domain, but using the modified
$\theta$ and $\phi$ variables.

We shall denote by $\theta = (\theta_{e_1}, \ldots, \theta_{e_r})$
the vector variable that puts together all the $\theta_e$
functions, \blue{for $e$ ranging through the set of all edges
$\{e_1, \ldots, e_r\}$. Similarly, let $\{\phi_1, \ldots,
\phi_p\}$ list all the functions $\phi_{P, e}$, for all vertices
$P$ and all edges $e$ containing $P$ we shall denote by $\phi =
(\phi_{1}, \ldots, \phi_{p})$} the vector variable that puts
together all the $\phi_{P, e}$ functions. We then introduce the
space $W^{k, \infty}(\Sigma\PP)$ as the space of functions $u :
\PP \to \CC$ of the form
\begin{multline*}
    u(x, y, z) = f(x, y, z, \theta, \phi) =
    f(x, y, z, \theta_{e_1}, \ldots, \theta_{e_r},
    \phi_{1}, \ldots, \phi_{p}), \\
    f \in W^{k, \infty}(\PP \times (0, 2\pi)^r \times
    (0, \pi)^{p}).
\end{multline*}
Thus $f$ above has $k$ bounded weak derivatives. We let $\CIS :=
\bigcap_k W^{k, \infty}(\Sigma\PP)$. The point of this definition
is that, for example, $\theta_e$ is a smooth function on $\PP$
that is not in $W^{k, \infty}(\PP)$ for $k \ge 1$. On the other
hand $\theta_e \in \CIS$, by definition.

One can show as in \cite{AIN, BNZ-ND} that there exists a
canonical Riemannian manifold $\Sigma(\PP)$ such that
$\CI(\Sigma(\PP)) = \CIS$, so our notation is justified. The
construction of a space with this property is not very intuitive.
However, at this point, we do not assign any significance to
$\Sigma\PP$, which should be regarded in this paper just as a
symbol. \blue{(Let us mention however, that, had we used curved
boundaries, then the desingularizations $\Sigma D_j$ of the faces
would have been necessary. See \cite{MazzucatoNistor2}.)}

\subsection{Vector fields and $\CIS$} We now establish several
technical properties of the functions in $\CIS$, especially in
relation to the vector fields (differentials) $r_{\PP}\pa_x$,
$r_{\PP}\pa_y$, and $r_{\PP}\pa_z$.

Let us notice first that it follows right away from the definition
that $\CIS$ is closed under addition and multiplication.

\begin{lemma}\label{lemma.rho.r}\
Let $P$ be a vertex of $\PP$, then $\rho_P \in \CIS$. Similarly,
let $e = [AB]$ be the edge of $\PP$ joining the vertices $A$ and
$B$, then $\tilde r_{e} := \rho_{A}^{-1} \rho_{B}^{-1} r_e \in
\CIS$. In particular, $r_\PP := \prod_e \tilde r_e \times \prod_P
\rho_P \in \CIS$.
\end{lemma}

This is proved using polar coordinates. Assume $P$ belongs to the
edge $e$, then $\rho_P = (\sin \phi_{P, e} \cos \theta_e)^{-1}x$,
where this is defined ($x$ stands for the first component
variable). Similar formulas for $\rho_P$ in terms of $y$ and $z$
then combine, using a partition of unity on $\RR^3 \smallsetminus
\{P\}$ with functions in $\CIS$, to define $\rho_P$ globally as an
element in $\CIS$.

Similarly, $\tilde r_e = \rho_A \sin \phi_{A, e}$, so
$\tilde r_e/ \rho_A$ is ``smooth'' near $A$. The same argument,
together with a partition of unity, shows that $r_e \in \CIS$. Our
result then follows from the fact that $\CIS$ is closed under
products, by definition.

\begin{lemma}\label{lemma.dist}\
Let $\vartheta(p)$ be the distance from $p$ to the union of the
edges of $\PP$. Then there exists $C > 0$ such that
$C^{-1}\vartheta(p) \le r_{\PP} \le C\vartheta(p)$ for all $p \in
\PP$.
\end{lemma}

This lemma is proved using the homogeneity properties of the functions
$\vartheta$ and $r_{\PP}$ close to the vertices and edges of
$\PP$. Using a compactness argument, it is enough to prove that the
ratio $r_{\PP}/\vartheta$ is bounded and bounded away from zero in the
neighborhood of each point. This allows us to assume that $\PP$ is
either a dihedral angle or an infinite cone. If $\PP$ is the dihedral
angle $0 < \theta < \alpha$, with $\alpha$ fixed, then
$r_{\PP}/\vartheta = 1$. If $\PP$ is a cone with center the origin,
let $\alpha_t$ be the dilation with center the origin and ratio
$t$. Then $r_{\PP}(\alpha_t(p)) = t\,r_{\PP}(p)$ and
$\vartheta(\alpha_t(p)) = t \,\vartheta(p)$.  This shows that the
ratio $r_{\PP}(p)/\vartheta(p)$ depends only on $p/|p|$.  Furthermore,
$r_{\PP}$ is a continuous function on the compact set
$\overline{\PP}\cap S^{n-1}$, and the lemma follows from this.

\begin{lemma}\label{lemma.deriv}\
We have that the functions \ $r_{e} \pa_x \theta_e$,\, $r_{e} \pa_y
\theta_e$,\, $r_{e} \pa_z \theta_e$,\, $\rho_P \pa_x \phi_{P, e}$,\,
$\rho_P \pa_y \phi_{P, e}$,\, $\rho_P \pa_z \phi_{P, e}$,\, $\pa_x
r_e$,\, $\pa_y r_e$,\, $\pa_z r_e$,\, $\pa_x \rho_P$,\, $\pa_y
\rho_P$,\, and\, $\pa_z\rho_P$ are all in $\CIS$.
\end{lemma}

To prove this, let us notice first that we can use any linear
system of coordinates $(x, y, z)$. In particular, for each of the
above calculations, we can assume that our cylindrical or
spherical coordinate system is aligned to the coordinate system
$(x, y, z)$. Then the result is simply an exercise in the
calculation of the partial derivatives of the cylindrical
coordinates $\theta$ and $r$ and of the spherical coordinates
$\phi$ and $\rho$.

\begin{corollary}\label{cor.rPP}\ We have \ $\pa_x r_{\PP},\,
\pa_y r_{\PP},\, \pa_z r_{\PP} \in \CIS$.
\end{corollary}

\begin{proof}\
Let us concentrate on $\pa_x$. We use the product rule to compute
the derivative of $r_{\PP}$. A summand containing
$\pa_x \rho_P$
is in $\CIS$ by Lemma \ref{lemma.deriv}. Let $e = [AB]$.
The other products are obtained by replacing $\tilde r_e :=
\rho_A^{-1} \rho_B^{-1} r_e$ with
\begin{equation*}
    \pa_x(\tilde r_e) = \rho_A^{-1} \rho_B^{-1}\pa_x(r_e)
    - \rho_A^{-1} \pa_x(\rho_A) \tilde r_e -
   \rho_B^{-1} \pa_x(\rho_B) \tilde r_e.
\end{equation*}
The factors of $\rho_A^{-1}$ and $\rho_B^{-1}$ then cancel out in
the product defining $r_{\PP}$ and all the remaining factors are
in $\CIS$ by Lemma \ref{lemma.deriv}.
\end{proof}

\begin{proposition}\label{prop.closed}\ If $u \in \CIS$, then
the functions $r_{\PP} \pa_x u$, $r_{\PP} \pa_y u$, and $r_{\PP}
\pa_z u$ are in $\CIS$.
\end{proposition}

\begin{proof}\ This follows from Lemma \ref{lemma.deriv} and
$r_{\PP} \in r_{e} \CIS \cap \rho_P \CIS$.
\end{proof}

Let us denote by $\DiffO{m}(\PP)$ the differential operators of order
$m$ on $\PP$ linearly generated by differential operators of the form
\begin{equation*}
    u (r_{\PP} \pa)^{\alpha} := u (r_{\PP} \pa_x)^{\alpha_1}
    (r_{\PP} \pa_y)^{\alpha_2} (r_{\PP} \pa_z)^{\alpha_3}, \quad
    |\alpha| := \alpha_1 + \alpha_2+ \alpha_3 \le m, \; u\in \CIS.
\end{equation*}
We agree that $\DiffO{m}(\PP) := \CIS$ and we shall denote
$\DiffO{\infty}(\PP) := \bigcup_{m}\DiffO{m}(\PP)$.  \blue{In case
of edges (no vertices), similar algebras were considered also by
Mazzeo, \cite{Mazzeo1, Mazzeo2}. Algebras more closely related to
ours appear in \cite{MelroseScattering}.} To get more insight into
the structure of $\DiffO{\infty}(\PP)$, we shall need two simple
calculations that we formalize in the following lemma, whose proof
is based on the fact that $\pa_j r_{\PP} \in \CIS$.

\begin{lemma}\label{lemma.2.calc}\
Let $\lambda \in \RR$ and let $\pa_j$ and $\pa_k$ stand for either
of $\pa_x$, $\pa_y$, or $\pa_z$. Then $r_{\PP}^{-\lambda} (r_{\PP}
\pa_j) r_{\PP}^{\lambda} - r_{\PP} \pa_j = \lambda \pa_j (r_{\PP})
\in \CIS,$ and
\begin{equation*}
    [r_\PP \pa_j, r_\PP \pa_k] := (r_\PP \pa_j) (r_\PP \pa_k)
    -(r_\PP \pa_k) (r_\PP \pa_j) = \pa_j (r_\PP) r_\PP \pa_k -
    \pa_k (r_\PP) r_\PP \pa_j \in \DiffO{1}(\PP).
\end{equation*}
\end{lemma}

Then we have the following simple but basic result.

\begin{proposition}\label{prop.alg}\
We have $\DiffO{k}(\PP)\DiffO{m}(\PP) \subset \DiffO{k+m}(\PP)$
and hence $\DiffO{\infty}(\PP)$ is an algebra.
\end{proposition}

\begin{proof}\ We shall prove by induction on $k + m$ that
$\DiffO{k}(\PP)\DiffO{m}(\PP) \subset \DiffO{k+m}(\PP)$.  Indeed, if
$k+m = 0$, then $k=m=0$ and the statement is clearly true because
$\CIS$ is closed under products.  Let us assume then that $k + m >
0$. We need to show that $u (r_{\PP} \pa)^{\alpha} v (r_{\PP}
\pa)^{\beta} \in \DiffO{k+m}(\PP)$ if $u, v \in \CIS$ and $|\alpha| :=
\alpha_1 + \alpha_2 + \alpha_3 =k$, $|\beta| := \beta_1 + \beta_2 +
\beta_3 = m$, where $\alpha = (\alpha_1, \alpha_2, \alpha_3)$
and $\beta = (\beta_1, \beta_2, \beta_3)$.

If $m =0$, then the relation
\begin{equation*}
    u (r_{\PP} \pa)^{\alpha} v = \sum u (r_{\PP} \pa)^{\alpha'}
    \big[r_{\PP} \pa_j (v) \big](r_{\PP} \pa)^{\alpha''}
\end{equation*}
for suitable $\alpha', \alpha''$ with $|\alpha'| + |\alpha''| =
k-1$, together with the induction hypothesis and with Proposition
\ref{prop.closed}, shows that $u (r_{\PP} \pa)^{\alpha} v \in
{\DiffO{k}(\PP)}$.

Let now $m$ be arbitrary. We shall proceed by a second induction
on $m$. The same argument as in the paragraph above allows us to
assume that $v = 1$. We can also assume that the monomial
$(r_{\PP} \pa)^{\alpha} (r_{\PP} \pa)^{\beta}$ is already ordered
in the standard way. Then, using Lemma \ref{lemma.2.calc}, we
commute $r_\PP \pa_j$, the last derivative in $(r_{\PP}
\pa)^{\alpha}$, with $r_\PP \pa_k$, the first derivative in
$(r_{\PP} \pa)^{\beta}$. Induction on $k + m$ for the terms
containing $\pa_j (r_\PP) r_\PP \pa_k$ and $\pa_k (r_\PP) r_\PP
\pa_j$ and induction on $m$ or $k + m$ for the term containing
$(r_\PP \pa_k) (r_\PP \pa_j)$ then complete the proof of the fact
that $\DiffO{k}(\PP)\DiffO{m}(\PP) \subset \DiffO{k+m}(\PP)$.
\end{proof}

The above proposition gives the following useful corollary.

\begin{corollary}\label{cor.order.m}\
If $P$ is a differential operator of order $m$ with smooth
coefficients, then $r_\PP^m P \in \DiffO{m}(\PP)$.
\end{corollary}

\begin{proof}\
It is enough to show that $r_\PP^m \pa^\alpha \in \DiffO{m}(\PP)$
if $\pa^\alpha = \pa_x^{\alpha_1} \pa_y^{\alpha_2}
\pa_z^{\alpha_3}$ with $|\alpha| = m$. We shall again proceed by
induction on $m$. The case $m =1$ is obvious. Let $\pa_j$ be the
first derivative in $\pa^\alpha$, so that $\pa^\alpha = \pa_j
\pa^{\alpha'}$. Then Lemma \ref{lemma.2.calc} and Corollary
\ref{cor.rPP} give
\begin{equation*}
    r_{\PP}^m \pa^\alpha - (r_{\PP} \pa_j) (r_{\PP}^{m-1}
    \pa^{\alpha'}) =  - (m-1) \pa_j(r_\PP) (r_{\PP}^{m-1}
    \pa^{\alpha'}) \in \DiffO{m-1}(\PP).
\end{equation*}
Then Proposition \ref{prop.alg} shows that $\DiffO{1}(\PP)
\DiffO{m-1}(\PP) \subset \DiffO{m}(\PP)$. This and the induction
hypothesis allows us to complete the proof.
\end{proof}

The proof of the above corollary also shows that
\begin{equation}\label{eq.rd}
    r_{\PP}^m \pa_x^{\alpha_1} \pa_y^{\alpha_2} \pa_z^{\alpha_3} -
    (r_{\PP} \pa_x)^{\alpha_1} (r_{\PP} \pa_y)^{\alpha_2} (r_{\PP}
    \pa_z)^{\alpha_3} \in \DiffO{m-1}(\PP), \quad |\alpha| = m.
\end{equation}

\section{Function spaces on $\PP$\label{sec.SS}}

We now recall and study the \BK\ spaces $\Kond{m}{a}(\PP) :=
\Kond{m}{a}(\PP; \vartheta)$ and $\Kond{m}{a}(\pa \PP; \vartheta)$ on
a $3$-simensional polyhedral domain $\PP$ and its boundary~$\pa\PP$.
These spaces are weighted Sobolev spaces with weight given by
$\vartheta$, the distance to the set of edges of $\PP$, as in Equation
\eqref{eq.def.WS}. Note that we can replace $\vartheta$ with $r_\PP$,
by Lemma \ref{lemma.dist} (we shall use this below).

\subsection{The \BK\ spaces}
We let
\begin{equation}
    W^{k, p, a}_{BK}(\PP) = \{u : \PP \to \CC,\
    r_{\PP}^{|\alpha| -a} \pa^{\alpha} u \in L^p(\PP), \text{
    for all } |\alpha| \le k\},
\end{equation}
for  $k \in \ZZ_+,\; a\in \RR,\; p\in [1,\infty]$.
If $p = 2$, we denote $\Kond{k}{a}(\PP) := W^{k, 2, a}_{BK}(\PP)$,
which coincides with the definition in the Introduction (Equation
\ref{eq.def.WS}).

We similarly define
\begin{multline*}
    W^{m, p, a}_{BK}(\pa \PP) = \{u : \pa \PP \to \CC,\
    r_{\PP}^{k-a} P (u\vert_{D_j}) \in L^p(D_j),\quad
    \text{for all } k \le m \\
    \text{ and all differential
    operators } P \text{ of order } k \text{ on } D_j,
    \; k \le m\}, \quad m \in \ZZ_+.
\end{multline*}
We let $\Kond{k}{a}(\pa \PP; \vartheta) := W^{m, 2, a}_{BK}(\pa
\PP)$. Thus $\Kond{k}{a}(\pa \PP; \vartheta) \simeq \bigoplus \Kond{k}{a}(D_j,
\vartheta)$ is thus a direct sum of weighted Sobolev spaces. Note that
we require {\em no compatibility conditions} for the resulting
functions on the faces $D_j$.

Equation~\eqref{eq.rd} and Lemma~\ref{lemma.dist} then give
immediately the following lemma.

\begin{lemma}\label{lemma.second.D}\
We have $\Kond{m}{a}(\PP) = \{u,\ \vartheta^{-a} Pu \in L^2(\PP),
\text{ for all } P \in \DiffO{k}(\PP)\}.$ A similar result holds for
$\Kond{m}{a}(\pa \PP; \vartheta)$ and for $W_{BK}^{k, p, a}(\PP)$.
\end{lemma}

Next, Proposition \ref{prop.closed} and Corollary
\ref{cor.rPP}, together with a straightforward calculation, show
the following.

\begin{lemma}\label{lemma.isom}\
The multiplication map $W^{m, \infty, b}_{BK}\times
\Kond{m}{a}(\PP) \ni (u, f) \mapsto uf \in \Kond{m}{a+b}(\PP)$ is
continuous. We also have $\CIS \subset W_{BK}^{m, \infty, 0}(\PP)$
and $r_{\PP}^b \in W_{BK}^{m, \infty, b}(\PP)$, and hence the map
$\Kond{m}{a}(\PP) \ni u \mapsto r_{\PP}^{b} u \in \Kond{m}{a+b}(\PP)$
is a continuous isomorphism of Banach spaces.
\end{lemma}

{}From this lemma we obtain right away the following result.

\begin{proposition}\label{prop.cont}\
Let $k \ge m$. Each $P_0 \in \DiffO{m}(\PP)$ defines a continuous
map $P_0 : \Kond{k}{a}(\PP) \to \Kond{k-m}{a}(\PP)$. The family
$r_{\PP}^{-\lambda} P_0 r_{\PP}^{\lambda}$ is a family of bounded
operators $\Kond{k}{a}(\PP) \to \Kond{k-m}{a}(\PP)$ depending
continuously on $\lambda$.

Similarly, if $P$ is a differential operator with smooth coefficients
on $\PP$, then $r_{\PP}^{-\lambda} P r_{\PP}^{\lambda}$ defines a
continuous family of bounded operators $\Kond{k}{a}(\PP) \to
\Kond{k-m}{a-m}(\PP)$.
\end{proposition}

\begin{proof} The first part follows from Lemma
\ref{lemma.second.D}. The second part follows from the first part of
this proposition and Lemma~\ref{lemma.isom}.
\end{proof}

We define the spaces $\Kond{-k}{a}(\PP)$, $k \in \ZZ_+$, by
duality. More precisely, let $\stackrel{\circ}{\Kond{k}{a}}(\PP)$
be the closure of $\CIc(\PP)$ in $\Kond{k}{a}(\PP)$. Then we
define $\Kond{-k}{-a}(\PP)$ to be the dual of
$\stackrel{\circ}{\Kond{k}{a}}(\PP)$, the duality pairing being an
extension of the bilinear form $(u, v) \mapsto \int_{\PP} uv\,
\dvol$. With this definition, we can drop the requirement that $k
\ge m$ in Proposition \ref{prop.cont}.

Let us also note that the resulting weighted Sobolev spaces on the
polygons $D_j$ are {\em different} from the weighted Sobolev spaces
obtained by using the distance to the vertices of these polygons. A
regularity theorem on $D_j$ would involve the latter weight (as in
Kondratiev's theorem \cite{Kondratiev67} mentioned in the
Introduction). A consequence of this is that the spaces
$\Kond{k}{a}(\pa \PP; \vartheta)$ behave more like the Sobolev spaces
defined on a smooth manifold without boundary than like the Sobolev
spaces defined on a bounded domain with (smooth) boundary. In
particular, we define $\Kond{-k}{-a}(\pa \PP; \vartheta)$ as the dual
of $\Kond{k}{a}(\pa \PP; \vartheta)$. The spaces $\Kond{s}{a}(\pa
\PP)$, $s \not\in \ZZ$, can be defined by interpolation, although in
this paper we shall use a different definition using partitions of
unity (see the following subsection; the two definitions are
equivalent, although we shall not need a proof of this fact in this
article).

\subsection{Definition of Sobolev spaces using partitions of unity}
As in \cite{AIN}, it is important to define the spaces
$\Kond{a}{m}(\PP)$ using partitions of unity. Similar
constructions were used in \cite{CGT, Shubin, Skrz, Triebel}. This
construction is possible because the spaces $\Kond{m}{3/2}(\PP)$
are the Sobolev spaces associated to the metric $r_{\PP}^{-2}g_E$,
where $g_E$ is the usual Euclidean metric.

We shall need the following lemma. Recall that $\vartheta(p)$
denotes the distance from~$p$ to the edges of $\PP$. In view
of Lemma \ref{lemma.dist}, in all estimates involving $\vartheta$,
we can replace $\vartheta$ with $r_{\PP}$, although not the other
way around, because $\vartheta$ is not smooth.

Let $\pa_{\operatorname{sing}} \PP$ be the union of the
edges of $\PP$ and $\PP' := \overline{\PP} \smallsetminus
\pa_{\operatorname{sing}} \PP$.

\begin{lemma}\label{lemma.Shubin}\
There is $\epsilon_0 \in (0, 1)$, an integer $\kappa$, and a sequence
$C_m > 0$ of constants such that, for any $\epsilon \in (0,
\epsilon_0]$, there is a sequence of points $\{x_j\} \subset \PP' :=
\overline{\PP} \smallsetminus \pa_{\operatorname{sing}} \PP$ and a
partition of unity $\phi_j \in \CIc(\PP')$ with the following
properties:
\begin{enumerate}[\rm (i)]
\item either $B(x_j, \epsilon \vartheta(x_j)/4)$ is contained in $\PP$
{\em or} $x_j \in \pa \PP$, $\vartheta(x_j) > 0$, and the ball $B(x_j,
\epsilon \vartheta(x_j))$ intersects only the face $D_i$ to which
$x_j$ belongs;
\item $\supp(\phi_j) \subset B(x_j, \epsilon \vartheta(x_j)/2)$ if
$x_j \in \pa \PP$ and  $\supp(\phi_j) \subset B(x_j, \epsilon
\vartheta(x_j)/8)$ otherwise;
\item $\phi_j(x_j) = 1$ and
$\|(r_\PP \pa)^{\alpha} \phi_j\|_{L^\infty(\PP)}\le
C_{|\alpha|}\epsilon^{-|\alpha|}$; and
\item a point $x \in \PP$ can belong to at most $\kappa$ of the
sets $B(x_j, \epsilon \vartheta(x_j))$.
\end{enumerate}
\end{lemma}

Let us notice that $\overline{B(x_j, \epsilon \vartheta(x_j))}$
does not intersect any edge of $\PP$ because $\epsilon < 1$.
Moreover, the conditions that $\|r_{\PP}\nabla \phi_j\|_{L^\infty} \le
C_1 \epsilon^{-1}$ and $\phi_j(x_j) = 1$ guarantee that the
support of $\phi_j$ is comparable in size with $\epsilon
\vartheta(x_j)$. This is reminiscent of the conditions appearing
in the definition of the Generalized Finite Element spaces
\cite{BabuBaOs, BabuCaOs, BabuMelenk1}.

A proof of this lemma will be given in the Appendix. It is
essentially a result that, in the case of non-compact manifolds,
goes back to Aubin. It was subsequently used by Gromov and in
\cite{AIN, Shubin, Skrz, Triebel}. We shall fix $\epsilon =
\epsilon_0$ in what follows and a sequence $x_j$ and a partition
of unity $\phi_j$ as in the lemma.

\begin{lemma}\label{lemma.f.sum}\
Let $u_k = \sum_{j=1}^k \phi_j u$, for $u \in
\Kond{m}{a}(\PP)$. Then $u_k \to u$ in $\Kond{m}{a}(\PP)$.
\end{lemma}

\begin{proof}\ Let $\Phi_k := \sum_{j = 1}^k\phi_j$.
We have that the sequence $(r_{\PP} \pa)^\alpha \Phi_k$ is bounded
in the `sup'-norm and converges to $0$ pointwise everywhere if
$\alpha \neq 0$. Similarly, $\Phi_k$ is bounded and converges to 1
pointwise everywhere. The result then follows from this using also
the Lebesgue dominated convergence theorem.
\end{proof}

Denote by $\alpha_j(x) = x_j + \vartheta(x_j)(x-x_j)$ be the
dilation of center $x_j$ and ratio $\vartheta(x_j)$. Let $J$ be
the set of indices $j$ such that $x_j \in \pa \PP$. Below,
by $H^m$ we shall mean either $H^{m}(\RR^3)$ or $H^{m}(\RR^3_+)$.
Also, denote by
\begin{multline}\label{eq.def.nu}
    \nu_{m, a}(u)^2 := \sum_{j} \vartheta(x_j)^{3-2a}
    \|(\phi_j u) \circ \alpha_j\|_{H^{m}}^2 \\
    := \sum_{j \not \in J} \vartheta(x_j)^{3-2a}
    \|(\phi_j u) \circ \alpha_j\|_{H^{m}(\RR^3)}^2 +
    \sum_{j \in J} \vartheta(x_j)^{3-2a}
    \|(\phi_j u) \circ \alpha_j\|_{H^{m}(\RR^3_+)}^2.
\end{multline}
We agree that $\|(\phi_j u) \circ \alpha_j\|_{H^{m}} = \infty$ if
$(\phi_j u) \circ \alpha_j \not \in H^{m}(\RR^3)$ (or if $(\phi_j
u) \circ \alpha_j \not \in H^{m}(\RR^3_+)$, respectively). Note
that the functions $(\phi_j u) \circ \alpha_j$ will all have
support contained in a fixed ball, namely, the ball $B(0,
\epsilon_0/2)$ of radius $\epsilon_0/2$ and center the origin.
Moreover, all derivatives $\pa^\alpha (\phi_j \circ \alpha_j)$ are
bounded for each fixed $\alpha$ and arbitrary $j$ by Lemma
\ref{lemma.Shubin}.

\begin{proposition}\label{prop.alt.desc}\
We have $u \in \Kond{m}{a}(\PP)$, $m \in \ZZ$, if, and only if,
$\nu_{m, a}(u) < \infty$. Moreover, $\nu_{m, a}(u)$ defines an
equivalent norm on $\Kond{m}{a}(\PP)$.
\end{proposition}

The proof of this Proposition is standard (see \cite[Lemma
2.4]{BabuNistor1}, \cite{AIN}, or \cite{Triebel}); for $m < 0$ one
also has to check that both definitions are compatible with
duality. We include a brief sketch below.

\begin{proof}\ Let us also introduce
\begin{equation*}%\label{eq.def.tnu}
    \tilde{\nu}_{m, a}(u)^2 := \sum_{j}
    \|(\phi_j u)\|_{\Kond{m}{a}(\PP)}^2.
\end{equation*}
Then the fact that $\vartheta(x)/\vartheta(x_j)$ and
$\vartheta(x_j)/\vartheta(x)$  are bounded by $(1-\epsilon)^{-1}$
on the ball $B(x_j, \epsilon \vartheta(x_j))$, for $\epsilon \in
(0, 1)$ and a change of variables shows that $\tilde {\nu}_{m, a}$
and $\nu_{m, a}$ define equivalent norms. It is then enough to
prove that $\tilde{\nu}_{m, a}(u)$ defines an equivalent norm on
$\Kond{m}{a}(\PP)$. For $m = 0$, this follows from the
inequalities
\begin{equation*}
    \|r_{\PP}^{-a} u \|^2_{L^2(\PP)} \le
    \kappa \tilde{\nu}_{0, a}^2 \le \kappa
    \|r_{\PP}^{-a} u \|^2_{L^2(\PP)}.
\end{equation*}
For arbitrary $m$, we use induction on $m$ and the fact that
$\sum_j |(r_{\PP}\pa)^\alpha \phi_j(p)|$ is bounded uniformly in
$p \in \PP$ for all $\alpha$.
\end{proof}

We proceed in the same way to study the spaces $\Kond{s}{a}(\pa \PP;
\vartheta)$, $s \in \RR$. Let us identify the plane containing each
face $D_k$ of $\PP$ with a copy of $\RR^2$. Then let
\begin{equation}\label{eq.def.nu2}
    \mu_{s, a}(u)^2 := \sum_{j \in J} \vartheta(x_j)^{2-2a}
    \|(\phi_j u) \circ \alpha_j\|_{H^{s}(\RR^2)}^2, \quad s \in
    \RR_+.
\end{equation}

Note that only the indices $j$ for which $x_j \in \pa \PP$ are
used above. Also, note that the power of $\vartheta{x_j}$ was
changed from $3-2a$ to $2-2a$.

Then we have an analogous description of the spaces
$\Kond{s}{a}(\pa \PP; \vartheta)$, $s \in \ZZ$.

\begin{proposition}\label{prop.alt.desc2}\
We have $u \in \Kond{s}{a}(\pa \PP; \vartheta)$ if, and only if,
$\mu_{s,a}(u) < \infty$. Moreover, $\mu_{s, a}(u)$ defines an
equivalent norm on $\Kond{s}{a}(\pa \PP; \vartheta)$, $s \in \ZZ$.
\end{proposition}

We can then define $\Kond{s}{a}(\pa \PP; \vartheta)$, $s \in \RR$, as the space
of functions $u$ for which $\mu_{s, a}(u) < \infty$ with the induced
norm. From this we obtain, by reducing to the Euclidean case, the
following Trace Theorem. Let $\pa_{\operatorname{sing}} \PP$ be the
union of the edges of $\PP$ and $\PP' := \overline{\PP} \smallsetminus
\pa_{\operatorname{sing}} \PP$, as above.

\begin{theorem}[Trace theorem]\label{thm.trace}\
The space $\CIc(\PP')$ is dense in $\Kond{m}{a}(\PP)$, $m \in
\ZZ_+$. The restriction to the boundary extends to a continuous,
surjective map $\Kond{m}{a}(\PP) \to \Kond{m-1/2}{a-1/2}(\pa \PP;
\vartheta)$ for $m \ge 1$. For $m = 1$, the kernel of this map is the
closure of $\CIc(\PP)$ in $\Kond{1}{a}(\PP)$.
\end{theorem}

\begin{proof}\
Clearly $\CIc(\PP') \subset \Kond{m}{a}(\PP)$, for any $m \in
\ZZ_+$ and any $a \in \RR$. To prove that it is a dense subspace,
let $u \in \Kond{m}{a}(\PP)$. By Lemma \ref{lemma.f.sum}, we may
assume that the support of $u$ does not intersect
$\pa_{\operatorname{sing}}\PP$ (replace $u$ with $u_k$ for some
$k$ large). Then we use the fact that $\CI(\overline{\Omega)}$ is
dense in $H^m(\Omega)$ for $\Omega$ a smooth, bounded domain and
the fact that the $H^m$-norm is equivalent to the norm on
$\Kond{m}{a}(\PP)$ when restricted to functions with support in a
fixed compact $K$ such that $K$ does not intersect any edge of
$\PP$ (\ie $K \cap \pa_{\operatorname{sing}} \PP = \emptyset$).

We have
\begin{multline*}
    \mu_{m-1/2, a-1/2}(u\vert_{\pa \PP})^2 := \sum_{j \in J}
    \vartheta(x_j)^{3-2a} \|(\phi_j u) \circ
    \alpha_j\vert_{\pa \PP}\|_{H^{m-1/2}(\RR^2)}^2 \\
    \le C \sum_{j \in J} \vartheta(x_j)^{3-2a}
    \|(\phi_j u) \circ \alpha_j\|_{H^{m}(\RR^3)}^2 \le
    C \nu_{m, a}(u),
\end{multline*}
and hence the restriction map $\Kond{m}{a}(\PP) \to
\Kond{m-1/2}{a-1/2}(\pa \PP; \vartheta)$ is defined and continuous for $m \ge
1$, by Propositions \ref{prop.alt.desc} and \ref{prop.alt.desc2}.
To prove that this map is continuous, let us fix a continuous
extension operator $E : H^{m-1/2}(\RR^2) \to H^{m}(\RR^3)$. By
rotation and translation, we extend this definition to an
extension operator $E : H^{m-1/2}(V) \to H^{m}(\RR^3)$, for any
two dimensional subspace $V \subset \RR^3$.

Let then $v : \pa \PP \to \CC$ be a function in
$\Kond{m-1/2}{a-1/2}(\pa \PP; \vartheta)$. Let us fix a function
$\psi \in \CIc(\RR^3)$ with support in the ball $B(0, \epsilon_0)$
of radius $\epsilon_0$ and center at the origin such that $\psi =
1$ on $B(0, \epsilon_0/2)$. Let $v_j(p) = \phi_j (\alpha_j(p)) u
(\alpha_j(p))$, which is defined on a subspace of $\RR^3$ of
dimension $2$. We define
\begin{equation*}
    u = \sum_j \big ( \psi E(v_j) \big)
    \circ \alpha_j^{-1}.
\end{equation*}
Then $u \in \Kond{m}{a}(\PP)$ and $u\vert_{\pa \PP} = v$.

Finally, let $u \in \Kond{1}{a}(\PP)$ such that $u\vert_{\pa \PP}
= 0$. Let $u_k$ be as in Lemma \ref{lemma.f.sum}. Then $u_k
\vert_{\pa \PP} = 0$. Using again the equivalence of the
$H^1$ and $\Kond{1}{a}(\PP)$--norms on functions with support in a
fixed compact set $K$ such that $K \cap \pa_{\operatorname{sing}}
\PP = \emptyset$, we see that each $u_k$ can be approximately in
$\Kond{1}{a}(\PP)$ as well as we want by a function $v_k \in
\CIc(\PP')$. Then we can take our approximation of $u$ to be $v =
\sum_{k=1}^N v_k$, for $N$ large enough.
\end{proof}

\section{Proof of the regularity theorem\label{sec.5}}

We include in this section the proof of Theorem~\ref{thm.princ}.
Its proof is reduced to the Euclidean case using a partition of
unity $\phi_j$ satisfying the conditions of Lemma
\ref{lemma.Shubin} for $\epsilon=\epsilon_0$, as in the previous
section.

\begin{proof} (of Theorem \ref{thm.princ}.)
The trace theorem, Theorem \ref{thm.trace} allows us to assume
that $u\vert_{\pa \PP} =0$. We then notice that, locally,
Theorem~\ref{thm.princ} is a well known statement. Namely, let us
consider a function $v$ with support in the ball of radius
$\epsilon_0$. We assume that either $v \in H^1(\RR^3)$ or $v \in
H_0^1(\RR^3_+)$ (that is, $v = 0$ on $\RR^2$, the boundary of
$\RR^3_+ = \{z \ge 0\}$). Then there exists a constant $C > 0$
such that, for all $m \ge 0$,
\begin{equation}\label{eq.est}
    \|v\|_{H^{m+1}(\RR^3)}^2 \le
    C_r \big(\|\Delta v\|_{H^{m-1}(\RR^3)}^2
    + \|v\|_{L^2(\RR^3)}^2 \big).
\end{equation}
or, respectively,
\begin{equation}\label{eq.est.bdry}
    \|v\|_{H^{m+1}(\RR^3_+)}^2 \le
    C_r \big(\|\Delta v\|_{H^{m-1}(\RR^3_+)}^2
    + \|v\|_{L^2(\RR^3_+)}^2 \big).
\end{equation}
The constant $C_r$ in the two equations above depends only on
$\epsilon_0$. (In fact, Equation \eqref{eq.est} implies
Equation \eqref{eq.est.bdry}, by taking $v$ to be odd with respect
to the reflection in the boundary of the half space $\RR^3_+$.)

We shall proceed by induction on $m \ge 0$. For $m = 0$, the
result is tautologically true, because of the term
$\|u\|_{\Kond{0}{a+1}(\PP)}$ on the right hand side of the
regularity estimate of Theorem \ref{thm.princ}. Let now
$\{\phi_j\}$ be the partition of unity and $\alpha_j$ be dilations
appearing in Equation \eqref{eq.def.nu}. In particular, the
partition of unity $\phi_j$ satisfies the conditions of Lemma
\ref{lemma.Shubin}, which implies that $\supp(\phi_j) \subset
B(x_j, \epsilon_0 \vartheta(x_j)/2)$ if $x_j \in \pa \PP$ and
$\supp(\phi_j) \subset B(x_j, \epsilon_0 \vartheta(x_j)/8)$
otherwise.  We also have that all derivatives of order $\le k$ of
the functions $\phi_j \circ \alpha_j$ are bounded. This implies in
turn that the commutator
\begin{equation*}
    P_j := [\Delta, \phi_j \circ \alpha_j]
    := \Delta (\phi_j \circ \alpha_j)
    - (\phi_j \circ \alpha_j) \Delta
\end{equation*}
is a differential operator all of whose coefficients have bounded
derivatives.

Let $\|v\|_{H^m}$ denote either $\|v\|_{H^m(\RR^3)}$ or
$\|v\|_{H^m(\RR^3_+)}$, depending on where the function $v$ is
defined. Let $\eta_j = \psi \circ \alpha_j^{-1}$, where $\psi \in
\CIc(\RR^3)$ has support in $B(0, \epsilon_0)$ and is equal to $1$
on $B(0, \epsilon_0/2)$, as before.

Then Equations \eqref{eq.est} and \eqref{eq.est.bdry} and the
above remarks give
\begin{multline*}
    \nu_{m + 2, a}(u)^2 :=  \sum_j \vartheta(x_j)^{3-2a}
    \|(\phi_j u) \circ \alpha_j\|_{H^{m+2}}^2 \\ \le C_r
    \sum_j \vartheta(x_j)^{3-2a} \Big( \|\Delta [(\phi_j u)
    \circ \alpha_j]\|_{H^{m}}^2 +
    \|(\phi_j u)\circ \alpha_j\|_{L^{2}}^2 \Big)
    \\ \le C \sum_j \vartheta(x_j)^{3-2a}
    \Big( \| (\phi_j \circ \alpha_j)
    \Delta (u \circ \alpha_j)\|_{H^{m}}^2
    + \|P_j (u \circ \alpha_j)\|_{H^{m}}^2 +
    \|(\phi_j u) \circ \alpha_j\|_{L^{2}}^2 \Big) \\ \le C
    \sum_j \vartheta(x_j)^{3-2a} \Big( \vartheta(x_j)^{4}\|
    (\phi_j \Delta u) \circ \alpha_j) \|_{H^{m}}^2
    + \|(\eta_j u) \circ \alpha_j\|_{H^{m+1}}^2
    + \|(\phi_j u) \circ \alpha_j\|_{L^{2}}^2 \Big) \\
    \le C \big( \nu_{m, a-2}(\Delta u)^2
    + \sum_j \nu_{m+1,a}(\eta_j u)^2
    + \nu_{0,a}(u)^2\big).
\end{multline*}
Since no more than $\kappa$ of the functions $\eta_j u$ are
non-zero at any given point of~$\PP$ and all the derivatives
$(r_\PP \pa)^\alpha \eta_j$ are bounded for all fixed $|\alpha|$,
we obtain that $\sum_j \nu_{m+1,a}(\eta_j u)^2 \le C
\nu_{m+1,a}(u)^2$. This then gives
\begin{equation*}
    \nu_{m + 2, a}(u)^2 \le C \big( \nu_{m, a-2}(\Delta u)^2 +
    \nu_{m+1,a}(u)^2\big).
\end{equation*}
By induction on $m$ we then obtain
\begin{equation*}
    \nu_{m + 2, a}(u)^2 \le C \big( \nu_{m, a-2}(\Delta u)^2 +
    \nu_{0,a}(u)^2\big).
\end{equation*}
The result then follows from Proposition~\ref{prop.alt.desc},
which states that the norms $\|\,\cdot\,\|_{\Kond{t}{a}(\PP)}$ and
$\nu_{t, a}$ are equivalent.
\end{proof}

See \cite{BNZ-I} for applications of these results, especially of the
above theorem.

By contrast, it is known that in the framework of the usual
Sobolev spaces $H^m(\PP)$, the smoothness of the solution of
\eqref{eq.BVP} is limited \cite{Dauge, Grisvard1, Grisvard2,
JerisonKenig, MitreaTaylor}.

\appendix
\section{Additional constructions}

In this appendix we explain how to modify the constructions of the
functions $\theta_e$ and $\phi_{P,e}$ introduced in Section~\ref{sec2}
when $\PP$ is not convex and how to construct a partition of unity
satisfying the conditions of Lemma~\ref{lemma.Shubin}.

\subsection{The modified functions $\theta$, $\phi$, and $r_e$}
We continue to denote by $\rho_P(p)$ the distance from $p$ to the
vertex $P$. By a dilation, we can assume that each edge of $\PP$
has length at least 4.

Let us first modify the functions $\phi_{A, e}$. \blue{We can find
$\delta > 0$ small enough so that for any vertex $P$, the sets
$\phi_{P, e} > \pi - 2\delta$ do not intersect ($e$ ranges through
the set of edges containing $P$).} Let $e = [AB]$ and $\psi_1 :
[0, \pi] \to [0, 1 - \delta]$ be a smooth, non-decreasing function
such that $\psi_1(x) = x$ for $0 \le x \le \pi - 2\delta$ and
$\psi_1(x) = \blue{\pi - \delta}$ for $x \ge \pi - \delta$. Also,
let $\psi_2 : [0, \infty) \to [0, 1]$ be a smooth, non-increasing
function such that $\psi_2(x) = 1$ for $0 \le x \le 1$ and
$\psi_2(x) = 0$ for $2 \le x $. Then we replace $\phi_{A, e}$ with
$\psi_1(\phi_{A, e}) \psi_2(\rho_A)$. This modifies the function
$\phi_{A, e}$ to make it smooth everywhere except on
$\overline{e}$.

\blue{We now modify the functions $\theta_e$. They will be
modified in two ways. Let us fix an edge $e = [AB]$. To understand
these modifications, it is useful to think of the spherical domain
$\omega_A$ associated to the vertex $A$. The old function
$\theta_e$ served the purpose of {\em both} desingularizing
$\omega_A$ close to the vertex associated to $e$ and of providing
global coordinates on $\omega_A$ away from the vertices (toghether
with the functions $\phi_{A, e}$). These two purposes of the old
$\theta_e$ will be accomplished by two modified functions
$\theta$. Let $\psi_3 : [0, \pi] \to [0, 1]$ be a smooth,
non-increasing function such that $\psi_3(x) = 1$ for $x \in [0,
\alpha]$ and $\psi_3(x) = 0$ for $x \ge 2 \alpha$. We then
similarly modify $\theta_e$ by replacing it with $\psi_3(\phi_{A,
e})\psi_3(\phi_{B, e}) \blue{\psi_2(r_e)} \theta_e$.  This will
make $\theta_e$ defined and smooth everywhere in space except on
$\overline{e}$ (if $\gamma$ is large enough). The resulting
function $\theta_e$ serves the purpose of desingularizing
$\omega_A$ near the vertex corresponding to $A$. Let next $\psi_4
: [0, 2\pi] \to [0, 2\pi]$ be a smooth function such that
$\psi_4(t) = t$ for $t \in [2\epsilon, 1-2\epsilon]$ and $\psi_4
(t) = 0$ for $t \in [0, \epsilon] \cup [1-\epsilon, 1]$. The
second kind of functions $\theta_A$ will be obtained by
considering $\psi_4(\theta_e) \psi_4(2\phi_{A, e})
\psi_4(2\phi_{B, e})$ for $\epsilon > 0$ small enough and {\em
all} choices of faces passing through $e$. (To define the old
functions $\theta_e$, we first chose a plane through $e$ and
containing one of the faces of $\PP$. This plane was the plane
where $\theta_e = 0$. For the new functions $\theta_e$, we
consider {\em all} the planes through $e$ and containing one of
the faces of $\PP$.) These new functions will be smooth on
$\omega_A$ near its vertices, but provide global coordinates away
from the vertices.}

Finally, let $\psi_3$ and $e = [AB]$ be as in the above paragraph.
We then replace $r_e$ with $\psi_3(\phi_{A, e})\psi_3(\phi_{B,
e})r_e + (1 - \psi_3(\phi_{A, e})) \rho_A + (1 - \psi_3(\phi_{B,
e})) \rho_B$.

\blue{Let us notice that one can define $\Sigma\PP$ directly,
which would provide the definition of $\CIS$ as the space of
smooth functions on $\Sigma\PP$ \cite{AIN}. The advantage of the
approach in \cite{AIN} is that it makes no distinction between the
cases when $\PP$ is convex or non-convex. The approach in this
paper has the advantage that it is much simpler and more intuitive
in the convex case.}

\blue{Further intuition in the construction of the spaces $\CIS$
can be obtained from the paper \cite{CDMaxwell2}, page 254, by
Costabel and Dauge where various regions and subregions of a
polyhedral domain were analyzed. See also \cite{NicaiseApel, BCD,
LubumaNicaise95}.}

\subsection{The partition of unity}

Our partition of unity will depend on parameters $(a, b, c)$ that
will be specified below.

First of all, let us denote by $B(P, n)$ the open ball of center
$P$ and radius $2^{-n} a$. By choosing $a$ small enough, we can
assume that the balls $B(P,1)$ do not intersect.  Then let $E_{e,
n}$ be the set of points $p \in \PP$ that do not belong to any
$B(P,2)$ and are at distance $\le 2^{-n}ab$ to the edge $e$. By
choosing $b$ small enough, we can assume that the sets $E_{e, 1}$
do not intersect. Let $\Omega_1$ be obtained from $\PP$ by
removing the sets $B(P,2)$ and $E_{e,2}$.

For each edge $e$, let $N_e$ be the plane normal to $e$. Project
$E_{e, 1}\smallsetminus E_{e, 2}$ onto $N_e$. The projection will
be the intersection of an annulus with an angle. Denote this
projection by $C_e$. We shall cover $C_e$ with disks of radius
$c/2$ and with disks of radius $c/8$. The disks of radius $c/2$
have the center on the straight sides of $C_e$ (the ones obtained
from the angle) and the disks of radius $c/8$ that have centers in
the interior of $C_e$ at distance at least $c/4$ to the angle
defining $C_e$. This yields the disks $D_1, \ldots, D_N$ with
centers $q_1,\ldots, q_N$.

Let $z$ be the variable along the line containing $e$. Then we cover
$E_{e, k+1} \smallsetminus E_{e, k+2}$, $k \in \ZZ_+$, with balls of
radius $2^{-k}c$ and centers of the form $(2^{-k}q_j, 2^{-k-3}c)$, if
$2^{-k}q_j$ is on one of the faces of $\PP$ and is inside $E_{e,
1}$. Otherwise, we consider the ball of radius $2^{-k-2}c$ with
centers of the form $(2^{-k}q_j, 2^{-k-3}c)$ as long as the center is
still inside $E_{e, 1}$.

Let us cover
\begin{equation*}
\Omega_1 := \PP \smallsetminus \Big( \bigcup_P B(P,2) \cup
\bigcup_e E_{e,2}\Big)
\end{equation*}
with finitely many balls of radius $c/2$ or radius $c/8$ with centers
in $\Omega_1$ such that the balls of radius $c/2$ have the centers on
the faces of $\PP$ and the balls of radius $c/8$ are at distance at
least $c/4$ to the faces of~$\PP$.

Let $D_{P,1}, \ldots, D_{P, N}, \ldots$ be the balls already
constructed with centers in $B(P, 1) \smallsetminus B(P, 2)$. Then
consider also the balls $2^{-k}D_{P,1}, \ldots, 2^{-k}D_{P, N},
\ldots$ obtained by dilations of ratio $2^{-k}$ and center $P$. We
repeat this construction for all vertices $P$ and all $k \in
\ZZ_+$. We consider all the balls $D_1, D_2, \ldots,$ constructed
so far (relabeled into a sequence) from the coverings of
$E_{e,1}$, $\Omega_1$, and from the dilations of ratio $2^{-k}$
for all the vertices $P$, as already explained. If we choose $c$
small enough (after the choices of $a$ and $b$ have been made as
explained above), then the sequence of these balls is locally
finite, the center of each ball is either on the faces of $\PP$ or
the closure of the ball is inside $\PP$. Moreover, for any such
ball $D$ with center $p$ and radius $r$, we have that
$r/\vartheta(p)$ is bounded from above and bounded from below from
zero, say $r/\vartheta(p) \in [\epsilon_0, \epsilon_0^{-1}]$, for
some $\epsilon_0 \in (0, 1)$. There is an integer $\kappa$ such
that no $\kappa+1$ of the balls constructed have a common point.

To any ball $D$ of center $q$ and radius $r$ we associate the {\em
bump function} $\psi_D(p) := \psi(|p - q|/r)$, where $\psi : [0,
\infty) \to [0,1]$ is smooth, is equal to 1 in a neighborhood of
$0$, is equal to $0$ in a neighborhood of $[1, \infty)$, and is
$>0$ on $[0, 1)$. Then we let $\eta = \sum \psi_{D_j}$ and $\phi_j
= \psi_{D_j}/\eta$. By further decreasing $c$, if necessary, we
see that our partition of unity (together with the points $x_j$
obtained as the centers of our balls) satisfies the conditions of
Lemma~\ref{lemma.Shubin} for the $\epsilon_0$ chosen above.

%\bibliographystyle{siam}
%\bibliographystyle{plain}
%\bibliography{num,babu}

\def\cprime{$'$} \def\cprime{$'$}
\def\ocirc#1{\ifmmode\setbox0=\hbox{$#1$}\dimen0=\ht0 \advance\dimen0
  by1pt\rlap{\hbox to\wd0{\hss\raise\dimen0
  \hbox{\hskip.2em$\scriptscriptstyle\circ$}\hss}}#1\else {\accent"17 #1}\fi}

\end{document}